\documentclass[10pt,oneside]{article}

\usepackage{amscd,amsfonts,amssymb,theorem}
\usepackage[centertags]{amsmath}
\usepackage[italian,english]{babel}
\usepackage[all]{xy}
\usepackage[dvips]{graphicx}

\usepackage[a4paper,hdivide={3.5cm,14cm, },vdivide={3cm,23.5cm, }]{geometry}

\newcommand{\Hodge}{\mathord{\mkern1mu *}}

\newcommand{\map}[3]{\mbox{${#1}\colon{#2}\to{#3}$}}

\newcommand{\R}{\ensuremath{\mathbb R}}
\newcommand{\C}{\ensuremath{\mathbb C}}
\newcommand{\Z}{\ensuremath{\mathbb Z}}

\newcommand{\Ha}{\ensuremath{\mathbb H}}
\newcommand{\Ca}{\ensuremath{\mathbb O}}
\newcommand{\Cp}[1]{\ensuremath{{\C\mathbb P}^{#1}}}

\newcommand{\iso}{\ensuremath{\simeq}}

\newcommand{\EndDim}{\ensuremath{\nopagebreak\hfill\blacksquare}}

\newcommand{\sform}{\ensuremath{\phi}}
\newcommand{\gform}{\ensuremath{\omega}}

\newcommand{\gtwo}{\ensuremath{\text{\upshape \rmfamily G}_2}}

\newcommand{\spin}[1]{\ensuremath{\text{\upshape \rmfamily Spin}(#1)}}
\newcommand{\gl}[1]{\ensuremath{\text{\upshape \rmfamily GL}(#1)}}

\newcommand{\so}[1]{\ensuremath{\text{\upshape \rmfamily SO}(#1)}}
\newcommand{\un}[1]{\ensuremath{\text{\upshape \rmfamily U}(#1)}}
\newcommand{\su}[1]{\ensuremath{\text{\upshape \rmfamily SU}(#1)}}
\newcommand{\symp}[1]{\ensuremath{\text{\upshape \rmfamily Sp}(#1)}}

\newcommand{\Hol}[1]{\ensuremath{\text{\upshape \rmfamily Hol}(#1)}}
\newcommand{\im}[1]{\ensuremath{\text{\upshape \rmfamily Im}(#1)}}
\newcommand{\isom}[1]{\ensuremath{\text{\upshape \rmfamily Isom}(#1)}}
\newcommand{\cone}[1]{\ensuremath{\mathcal{C}(#1)}}

\newtheorem{lem}{Lemma}[section]

\newtheorem{teo}[lem]{Theorem}
\newtheorem{cor}[lem]{Corollary}
\newtheorem{pro}[lem]{Proposition}
{\theorembodyfont{\rmfamily} \newtheorem{rem}[lem]{Remark}
            
            \newtheorem{exa}[lem]{Examples}
            }

\newenvironment{D}[1][]{{\nopagebreak\noindent\em Proof#1: }}{\EndDim}

\newcommand{\lee}{\ensuremath{\theta}}
\newcommand{\leespin}{\ensuremath{\Theta}}

\newcommand{\lie}{\ensuremath{\mathcal{L}}}

\date{}

\title{Locally conformal parallel \gtwo\ and \spin{7} manifolds}
\author{Stefan Ivanov, Maurizio Parton, Paolo Piccinni}

\def\sideremark#1{\ifvmode\leavevmode\fi\vadjust{\vbox to0pt{\vss
 \hbox to 0pt{\hskip\hsize\hskip1em
 \vbox{\hsize2.5cm\tiny\raggedright\pretolerance10000
 \noindent #1\hfill}\hss}\vbox to8pt{\vfil}\vss}}}%

\begin{document}

\maketitle

\begin{abstract}
We characterize compact locally conformal parallel \gtwo\ (respectively, \spin{7}) manifolds as fiber bundles over $S^1$ with compact nearly K\"ahler (respectively, compact nearly parallel \gtwo) fiber.
A more specific characterization is provided when the local parallel structures are flat.
\end{abstract}

\section{Introduction}

Locally conformal parallel structures have been studied for a long
time with respect to different natural geometries. The K\"ahler
condition is the oldest example of parallel structure and locally
conformal K\"ahler geometry is by now a quite developed subject
with several interactions to other geometries
\cite{DrOLCK,OrnLCK}. Other choices include locally conformal
hyperk\"ahler and quaternion K\"ahler metrics \cite{OrPLCK}. In
particular, the useful technique of the quotient construction,
coming from symplectic geometry, has been described both for
locally conformal K\"ahler and locally conformal hyperk\"ahler
metrics \cite{GOPLCK,GOPRVS}. By looking at further groups whose
holonomy has significance, the choices of \gtwo\ and \spin{7}\
appear as deserving of attention. 

In fact locally conformal
parallel \gtwo\ and \spin{7} structures are both natural classes
in the frameworks of the the symmetries considered in
%
%
%
\cite{FeGRMS,FerCRM}. More recently, conformal parallel \gtwo\
structures have been studied in details on solvmanifolds \cite{ChFCPG}.
A locally conformal parallel structure is encoded in the Lee form
\lee, the 1-form representing the irreducible component of the
covariant derivative of the fundamental form in the standard
representation. When \lee\ can be assumed to be parallel for one
metric in the conformal class the geometry of the manifold is
described, in the K\"ahler and hyperk\"ahler settings, by a
fibration over a circle \cite{OrVSTC,VerVTL,GOPRVS}. In fact, we obtain here
the corresponding characterization in the \gtwo\ and \spin{7} cases.

\medskip
\noindent{\bf Theorem A}
A compact Riemannian 7-manifold $M$ admits
a locally conformal parallel
\gtwo\ structure if and only if
there exists a fibre bundle
$M\rightarrow S^1$ with abstract fibre $N/\Gamma$, where $N$ is a
compact simply connected nearly K\"ahler 6-manifold and $\Gamma$
is a finite subgroup of complex isometries of $N$
acting freely. Moreover, the cone
$\cone{N/\Gamma}$ covers $M$ with cyclic infinite covering
transformations group.

\medskip
\noindent{\bf Theorem B}
A compact Riemannian 8-manifold $M$ admits
a locally conformal parallel
\spin{7} structure if and only if
there exists a fibre bundle
$M\rightarrow S^1$ with abstract fibre $N/\Gamma$, where $N$ is a
compact simply connected nearly parallel \gtwo\ manifold and $\Gamma$
is a finite subgroup of isometries of $N$ acting freely.
Moreover, the cone
$\cone{N/\Gamma}$ covers $M$ with cyclic infinite covering
transformations group.

\medskip
When the local parallel structures are flat then this statement
simplifies. In the \gtwo\ case, one gets a trivial circle bundle whose fibre
is the standard nearly K\"ahler 6-sphere $S^6$
({\bf Theorem~\ref{g2flat}}). In the \spin{7} case, the fibre is a
spherical space form $S^7/\Gamma$, where $\Gamma$ is a finite
subgroup of \spin{7} ({\bf Theorem~\ref{sflat}}).
In particular, the problem of describing all locally conformal parallel flat
\spin{7} structures is thus reduced to the possibility of recognizing, inside the list
of finite subgroup of \so{8}, the ones which are contained in \spin{7}.

An independent proof of Theorem A has recently been given by M.~Verbitsky in \cite{VerIVF}. Moreover, the statements given here concerning the universal coverings have been recently formulated for manifolds of any dimension in the setting of geometric structures of vectorial type by I.~Agricola and T.~Friedrich in \cite{AgFGSV}.

\bigskip
\par \noindent
{\bf Acknowledgements}. The first author thanks I.N.d.A.M.\ and
Universit\`a di Roma ``La Sapienza'' for hospitality and support.
All the authors thank Centro Ennio De Giorgi in Pisa for support and the warm atmosphere during the
program ``Differential Geometry and Topology'' in Fall 2004. The
authors are grateful to Kris Galicki, Rosa Gini, Liviu Ornea,
Thomas Friedrich, Simon Salamon, Misha Verbitsky for suggestions and
discussions.

\section{Preliminaries}

Let $e_1,\dots,e_7$ be an oriented
orthonormal basis of $\R^7$.
The subgroup of $\gl{7}$ fixing the 3-form
\begin{equation}
\gform =e_{127} - e_{236} + e_{347}+e_{567} - e_{146} - e_{245} + e_{135}\label{11}
\end{equation}
is the compact Lie
group $\gtwo\subset\so{7}$,
and $\gform$ corresponds to a real spinor so that
$\gtwo$ can be identified with the isotropy group of a non-trivial
real spinor. 

A Riemannian manifold $(M^7,g)$ is a $\gtwo$ manifold
if its structure group reduces to $\gtwo$, or equivalently if there exists a 3-form $\gform$ on $M$ locally given by \eqref{11}.
A $\gtwo$ manifold is said to be parallel if $\gform$ is
parallel (namely, if $\Hol{g}\subset\gtwo$) and \emph{locally conformal parallel}\/ if the 1-form \lee\ defined by 
\[
\lee=-\frac{1}{3}\Hodge(\Hodge d\gform\wedge\gform)
\]
satisfies
\begin{equation}\label{lcpg2}
d\gform=\frac{3}{4}\lee\wedge\gform,\qquad d\Hodge\gform=\lee\wedge\Hodge\gform.
\end{equation}
The latter terminology is justified by the fact that if \eqref{lcpg2} holds, then \lee\ is closed \cite{CabRMG} and the functions locally defined by $\lee=df$ gives rise to local parallel \gtwo\ structures $e^{3f}\gform$.
As usual, the choice of constants is a matter of habits and convenience.

Finally, a $\gtwo$ structure with co-closed Lee form is called a \emph{Gauduchon}\, $\gtwo$ structure.

Next, let
$e_0,\dots,e_7$ be an oriented orthonormal basis of ${\R}^8=\R\oplus\R^7$.
The subgroup of $\gl{8}$ fixing the 4-form
\begin{equation}
\sform = e_0\wedge\gform +\Hodge_7\gform\label{1}
\end{equation}
is the compact Lie
group $\spin{7}\subset\so{8}$.
Like in the $\gtwo$ case, $\sform$ corresponds to a real spinor and $\spin{7}$ can be identified
with the isotropy group of a non-trivial real spinor.

A $\spin{7}$ structure on $(M^8,g)$ is a
reduction of its structure group to $\spin{7}$, or equivalently a nowhere vanishing
4-form $\sform$ on $M$, locally written as \eqref{1}.
A $\spin{7}$ manifold $(M,g,\sform)$ is said to be parallel  if
$\Hol{g} \subset\spin{7}$, or equivalently if $\sform$ is parallel.
The Lee form of a \spin{7} manifold is now
\[
\leespin = -\frac{1}{7}\Hodge(\Hodge
d\sform\wedge\sform)
\]
and if $d\sform = \leespin\wedge\sform$ then $d\leespin=0$ \cite{CabRMS} and $\sform$ is locally conformal
parallel. Again, a $\spin{7}$ structure with co-closed Lee form will be called a
Gauduchon $\spin{7}$ structure.

\section{Structure of compact manifolds}

We will need the following fact concerning the
homotheties of a cone over a compact Riemannian manifold.
\begin{teo}\label{gopp} {\rm\bf\cite{GOPRVS}}
Let $N$ be a compact Riemannian manifold, and denote by \cone{N}
its Riemannian cone. Let $\Gamma$ be a discrete subgroup of homotheties
of \cone{N} acting freely and properly discontinously. Then $\Gamma$ is
a finite central extension of \Z, and the finite part is a finite
subgroup of \isom{N}.
\end{teo}

Recall that an almost Hermitian $(N^6,J,g)$ with K\"ahler form $F$
is \emph{nearly K\"ahler} if $dF$ decomposes as a
$(3,0)+(0,3)$-form. Similarly a $\gtwo$ manifold $(N^7,\gform,g)$
is nearly parallel if $d\gform=\lambda\Hodge\gform$ with constant
$\lambda>0$. In both cases $g$ is Einstein with positive scalar
curvature.  It is easy to see that $N$ is nearly K\"ahler or
nearly parallel $\gtwo$ if and only if its Riemannian cone
\cone{N} is a parallel \gtwo\ or \spin{7} manifold,
respectively \cite{BarRKS}.

\subsection{Proof of Theorem~A}
Denote by
$[\gform]$ the conformal class of the fundamental 3-form, and observe that we can
more appropriately define locally conformal parallel \gtwo\ manifolds as pairs
$(M,[\gform])$ such that any representative of the conformal class
is locally conformal parallel. Of course, one has an associated
conformal class $[g]$ of Riemannian metrics. Note also that any $\gtwo$ manifold $M$ can be viewed
through its Lee form as a Weyl manifold. Thus, whenever $M$ is
compact, a Gauduchon $\gtwo$ structure can be obtained
\cite{GauSWE,FrIKSE}.

\begin{exa}
Consider a nearly K\"ahler manifold $N^6$, with fundamental form $F$.
Then $N\times S^1$
admits a locally conformal parallel $\gtwo$ structure defined by \cite{ChSITS}
\begin{equation}\label{flat}
\gform=dt\wedge F+dF.
\end{equation}
Thus the four known compact nearly K\"ahler
6-manifolds, namely $S^6$, $S^3\times S^3$, $\Cp{3}$ and the flag
manifold $F=\un{3}/(\un{1}\times \un{1}\times \un{1})$ give rise
to products with $S^1$ that admit a locally conformal parallel $\gtwo$ structure
\cite{HitSFS,SalAPS}. Further examples can be obtained by choosing $N$ to be the twistor
orbifold of the weighted projective planes $\Cp{2}_{p_0,p_1,p_2}$, the case $(p_0,p_1,p_2)=(1,1,1)$
corresponding to the flag manifold $F$. These orbifolds $N$ admit
a nearly K\"ahler structure and therefore their products with
$S^1$ are locally conformal parallel $\gtwo$ manifolds.
\end{exa}

\begin{pro}\label{tec}
Let $(M,[\gform])$ be a compact locally conformal parallel $\gtwo$
manifold. If $\gform$ is a Gauduchon
$\gtwo$ form,  $\lee$ its Lee form and $g$ the metric
with Levi-Civita connection $\nabla^g$, then $\nabla^g\lee=0$.
Moreover:
\begin{enumerate}
\item The Riemannian universal cover of $M$ is isometric, up to
homotheties, to the product $N\times\R$, where $N$ is a compact
6-dimensional Einstein manifold with positive scalar curvature.
\item The scalar curvature of $g$  is the positive constant
$s=\frac{15}{8}|\lee|^2$.
\item  The Lee flow preserves the Gauduchon $\gtwo$ structure,
that is $\lie_{\lee^\sharp}\gform=0$, where $\lee^\sharp$ denotes
the Lee vector field. 
\item If $\nabla$ is the unique
$\gtwo$ connection with totally skew-symmetric torsion $T$, then
$\nabla T=0$.
\end{enumerate}
\end{pro}

\begin{D}
The Weyl structure of $M$ defines a Weyl connection $D$ such that
$D g=\lee \otimes g$. Since $\lee$ is closed and the Ricci tensor
vanishes, Theorem 3 of \cite{GauSWE} can be applied. This gives
$\nabla^g \lee = 0$ as well as property i). Statement ii) is a
direct consequence of the formula (4.19) in \cite{FrIKSE}.
To prove iii), observe that $\lie_{\lee^\sharp} \gform=\nabla_{\lee^\sharp}^g
\gform$. Using the fact that locally $\lee =-4df$ and that the Weyl
connection $D$ coincides locally with the Levi-Civita connection of
the parallel $\gtwo$ structure $(\gform'=e^{3f}\gform,g'=e^{2f}g)$
we get
\begin{equation}\label{wg2}
D_{\lee^\sharp}X=\nabla^g_{\lee^\sharp}X-\frac14|\lee|^2X, \qquad D \gform=\frac34\lee\otimes\gform.
\end{equation}
On the other hand, using the first equation in \eqref{wg2} we
calculate $D_{\lee^\sharp}\gform=\nabla^g_{\lee^\sharp}\gform +
\frac34|\lee|^2\gform$. Compare with the second equation in
\eqref{wg2} to conclude
$0=\nabla^g_{\lee^\sharp}\gform=\lie_{\lee^\sharp}\gform$.
Finally iv) follows from the fact that the unique $G_2$ connection
having as torsion a 3-form is given by \cite{FrIKSE}
\[
g(\nabla_XY,Z)=g(\nabla^g_XY,Z)
+\frac12T(X,Y,Z),
\]
where $T=\frac14\Hodge(\theta \wedge\gform
)=-\frac14i_{\theta}(\Hodge\gform)$. The last two equalities lead
to $\nabla\theta=0$ and hence to $\nabla T=0$.
\end{D}

\begin{cor}\label{covering}
A 7-dimensional compact Riemannian manifold
$M$ carries a locally conformal parallel
\gtwo\ structure if and only if $M$ admits a covering which is a Riemannian
cone over a compact nearly K\"ahler 6-manifold and such that the
covering transformations are homotheties preserving the corresponding parallel
\gtwo\ structure.
\end{cor}

\begin{D}
According to Proposition~\ref{tec}, the universal covering $\tilde
M$ is isometric to $N\times \R$ with the product metric
$\tilde g=g_N+dt^2$ and $\tilde\lee=dt$. The local conformal
transformations making the structures parallel give rise to a
global $G_2$ structure on $\tilde M$. Therefore
\begin{gather*}
\tilde g'=e^{2t}\tilde g=e^{2t}\left(g_N+dt^2\right), \qquad
\tilde\gform'=e^{3t}\tilde\gform,
\end{gather*}
defines a parallel \gtwo\ structure on $\tilde{M}$. Observe that 
$\tilde g'$ is the cone metric and use \cite{BarRKS} to get the conclusion.
\end{D}

\begin{rem}\label{topology}
The universal covering used in the proof of
Corollary~\ref{covering} can be replaced by any Riemannian
covering of $M$ such that the Lee form is exact,
say $\tilde{\lee}=df$. Then \map{f}{\tilde{M}}{\R} is
a Riemannian submersion, because $|\tilde{\lee}|$ is constant.
Moreover, $f$ is globally trivial, for the Lee flow acts by isometries
and $\tilde{M}$ is complete.
\end{rem}

The previous Remark together with Theorem~\ref{gopp} gives that
among all the cones covering a compact $M$, one can always choose
the one whose covering transformations group is cyclic infinite.
\begin{pro}\label{structureg2}
Let  $(M,[\gform])$ be a compact locally conformal parallel
\gtwo\ manifold. Then there exists a covering
$\tilde{M}\stackrel{\Z}{\rightarrow} M$, where $\tilde{M}$ is a
globally conformal parallel \gtwo\ manifold.
\end{pro}

\begin{D}
As a consequence of the property $\nabla^g \lee = 0$, we can
always assume the Lee form to be of constant length and harmonic.
Moreover, the Lee field preserves \gform, so that it is Killing
and the Lee flow acts as $G_2$ isometries. Next,
Corollary~\ref{covering} gives a globally conformal parallel
covering $\cone{N}\stackrel{\Gamma}{\rightarrow} M$, where $N$ is
a compact nearly K\"ahler 6-manifold. Denote by
\map{\rho}{\Gamma}{\R^+} the map given by the dilation factors of
elements of $\Gamma$. The isometries in $\Gamma$ are then
$\ker\rho$ and by Theorem~\ref{gopp} we obtain
\[
\Gamma/\ker\rho\iso\im\rho\iso\Z.
\]
Define $\tilde{M}\stackrel{\Z}{\rightarrow} M$ as
$\cone{N/\ker\rho}\stackrel{\Gamma/\ker\rho}{\rightarrow} M$.
\end{D}

\bigskip
If $\cone{N}\stackrel{\Z}{\rightarrow}M$ is the covering given by
Proposition~\ref{structureg2} and $\pi$ the projection of the cone
onto its radius, it is easy to check that $\pi$ is equivariant
with respect to the action of the covering maps on $\cone{N}$ and the
action of $n\in\Z$ on $t\in\R$ given by $n + t$.
From all of this, Theorem A follows.

\subsection{Proof of Theorem B}

Again, denote by
$[\sform]$ the conformal class of the fundamental 4-form, and note that locally conformal parallel \spin{7}
manifolds can be viewed as pairs
$(M,[\sform])$ such that any representative of the conformal class
is locally conformal parallel. Of course, one has an associated
conformal class $[g]$ of Riemannian metrics.

\begin{rem}\label{np1}
Like in the $\gtwo$ case, any compact $\spin{7}$ manifold admits a
Gauduchon $\spin{7}$ structure \cite{IvaCTP}. In fact, if
$(N^7,\gform)$ is any nearly parallel $\gtwo$ manifold, the
product $N^7\times S^1$ carries a natural locally conformal
parallel $\spin{7}$ structure defined by
\[
\sform=dt\wedge\gform + \Hodge_{N^7}\gform
\]
where $\Hodge_{N^7}$ is the Hodge star operator on $N^7$.
\end{rem}

\begin{exa}
Through the above formula one obtains several classes of examples of locally conformal parallel
\spin{7} manifolds. In fact, according to \cite{FKMNPG}
nearly parallel \gtwo\ manifolds split in 4 classes, in correspondance with 
the number of Killing spinors. These 4 classes are (I) proper
nearly parallel \gtwo\ manifolds, (II) Sasakian-Einstein manifolds,
(III) 3-Sasakian manifolds and (IV) the class containing only
$S^7$. In particular, the class (I) contains all squashed
3-Sasakian metrics, including $S_{\text{sq}}^7$ with its Einstein
metric. Moreover, we have in this class the homogeneous real
Stiefel manifold $\so{5}/\so{3}$ and the Aloff-Wallach spaces
\cite{CMSCGS}. The class (II) includes $S^7$
with any of its differentiable structures and any Einstein metric in the families constructed in \cite{BGKEMS}. The
class (III) is known to contain examples with arbitrarily large second Betti number
\cite{BGMCSM}.
Note that the topological $S^7$ can be endowed with nearly
parallel \gtwo\ structures in any of the above mentioned 4
classes.

These 4 classes correspond to the sequence of inclusions
$\spin{7}\supset \su{4}\supset \symp{2}\supset \{\text{id}\}$.
From this point of view, the products $N\times S^1$ carry a proper
locally conformal parallel \spin{7} structure, a locally conformal
K\"ahler Ricci-flat metric, a locally conformal hyperk\"ahler
metric and a locally conformal flat metric according to whether
$N$ belongs to the class (I), (II), (III) or (IV).

One can easily check that if $I,J,K$ denote the hypercomplex
structure on $N\times S^1$ (class (III)), then the \spin{7}
form $\sform$ is related to the K\"ahler forms
$\omega_I,\omega_J,\omega_K$ by $
2\sform=\omega_I^2+\omega_J^2-\omega_K^2$.
Thus the locally conformal hyperk\"ahler
property $d\omega_I=\lee\wedge\omega_I,\dots$ implies
$\lee=\leespin$. This holds of course also for $N$ in class
(II), where the conformality factors $f_\alpha$ on $N\times
S^1$ can be choosen in the same way for both the local K\"ahler
metrics and the local parallel \spin{7} forms.
\end{exa}

\begin{rem}
The existence of a Killing spinor on a nearly parallel $\gtwo$
manifold turns out to be equivalent to the existence of a parallel
spinor with respect to the unique $\gtwo$ connection with totally
skew-symmetric torsion \cite{FrIPSC}. The different classes of
structures correspond to the fact that the holonomy group of the
torsion connection is contained in $\symp{1}$, $\su{3}$ and
$\gtwo$, respectively.
\end{rem}
An easy consequence of the considerations in \cite{IvaCTP} is that
the unique $\spin{7}$ connection with totally skew-symmetric
torsion on a locally conformal parallel \spin{7} manifold is determined by the formula
\[
g(\nabla_XY,Z)=g(\nabla^g_XY,Z) + \frac12T(X,Y,Z),\quad
\text{where}\quad T =-\frac16 \Hodge(\leespin\wedge\sform)
=-\frac16i_{\leespin^\sharp}(\Hodge\sform).
\]
Now, in the exactly same way as in the \gtwo\ case, we have the
following.
\begin{pro}\label{tecs}
Let $(M^8,[\sform])$ be a compact locally conformal parallel
$\spin{7}$ manifold, where $\sform$ is a Gauduchon
$\spin{7}$ structure. Let $\leespin$ be its Lee form, $g$ be the
metric and $\nabla^g$ its Levi-Civita connection. Then
$\nabla^g\leespin=0$. Moreover:
\begin{enumerate}
\item The Riemannian universal cover is isometric, up to homothety,
to the product $N^7\times\R$, where $N$ is a compact 7-dimensional
Einstein manifold with positive scalar curvature.
\item The scalar curvature of $g$  is the positive constant
$s=\frac{21}{36}|\leespin|^2$.
\item The Lee flow preserves the Gauduchon $\spin{7}$ structure,
that is ${\lie}_{\leespin^\sharp}\sform=0$, where $\leespin^\sharp$ denotes
the Lee vector field.
\item If $\nabla$ denotes the unique $\spin{7}$ connection with
totally skew-symmetric torsion $T$, then $\nabla T=0$.
\end{enumerate}
\end{pro}

Moreover, similarly to the \gtwo\ case, we get:

\begin{cor}\label{coveringspin}
A 8-dimensional compact Riemannian manifold
$M$ carries a locally conformal parallel
\spin{7} structure if and only if $M$ admits a covering which is a Riemannian
cone over a compact nearly parallel \gtwo\ manifold and such that the
covering transformations are homotheties preserving the corresponding parallel
\spin{7} structure.
\end{cor}

\begin{pro}\label{structure}
Let  $(M,[\sform])$ be a compact locally conformal parallel
\spin{7} manifold.
Then there exists a covering $\tilde{M}\stackrel{\Z}{\rightarrow} M$,
where $\tilde{M}$ is
a globally conformal parallel \spin{7} manifold.
\end{pro}
Then Theorem~B follows.

\subsection{Locally conformal flat structures and further remarks}

The following statements specialize Theorems A and B whenever the local metrics are flat.

\begin{teo}\label{g2flat}
Let $(M,[\gform])$ be a compact locally conformal parallel flat
\gtwo\ manifold. Then $(M,[\gform])=(S^6\times
S^1,[\gform_0=\eta\wedge
F_0+dF_0])$, where $F_0$ is the
standard nearly K\"ahler structure on $S^6$.
\end{teo}

\begin{D}
The flatness of the local parallel \gtwo\ structures is equivalent
to the flatness of the Weyl connection of $M$, and from Corollary
at page 11 in \cite{GauSWE} the associated compact simply
connected nearly K\"ahler manifold is isometric to the round
sphere $S^6$. Thus Theorem A implies
that $M$ is a fibre bundle over $S^1$ with fibre $S^6/\Gamma$,
where $\Gamma$ is a finite subgroup of \so{7}. Since $M$ is
orientable, $\Gamma$ must be trivial, and the connectedness of
\so{7} implies that $M$ is isometric to the trivial bundle
$S^6\times S^1$. The \gtwo\ structure is given by
$[\gform]=[\eta\wedge F+dF]$, where $F$ is a nearly K\"ahler
structure on the round sphere $S^6$. But the standard
nearly K\"ahler structure on $S^6$ inherited from the imaginary
octonions is the only nearly
K\"ahler structure on $S^6$ compatible with the round metric \cite{FriPrC}, and this ends the
proof.
\end{D}


\begin{teo}\label{sflat}
Let $(M,[\sform])$ be a compact locally conformal parallel flat
\spin{7} manifold. Then $(M,[\sform])=(N\times
S^1,[\sform_0=\eta\wedge\gform_0+
\Hodge_{N}\gform_0])$, where $N=S^7/\Gamma$
is a spherical space form and $\gform_0$ is the
standard nearly parallel \gtwo\ structure on $S^7$.
\end{teo}

\begin{D}
As for the \gtwo\ case, the flatness of the local parallel
\spin{7} structures gives that the associated compact simply
connected nearly parallel \gtwo\ manifold is isometric to the
round sphere $S^7$. Thus Theorem B
implies that $M$ is a fibre bundle over $S^1$ with fibre
$S^7/\Gamma$, where $\Gamma$ is a finite subgroup of \so{8}, and the connectedness of \so{8} implies
that $M$ is isometric to the trivial bundle $(S^7/\Gamma)\times
S^1$. The \spin{7} structure is given by
$[\sform]=[\eta\wedge\gform + \Hodge_7\gform]$, where $\gform$ is
a $\Gamma$-invariant nearly parallel \gtwo\ structure on the round
sphere $S^7$ and $\Hodge_7$ is the associated Hodge star operator.
But the standard nearly parallel \gtwo\ structure on $S^7$
inherited from the standard \spin{7} structure on $\R^8$ is the
only nearly parallel \gtwo\ structure on $S^7$ compatible with the
round metric \cite{FriPrC}, and this completes the proof.
\end{D}

\begin{rem} Note that $\R^8-\{0\}$ admits a $\spin{7}$ structure that is not in the
conformal class generated by the standard one. This can be seen by looking at the proper nearly
parallel $\gtwo$ structure on the squashed sphere $S^7_\text{sq}$ and hence at
the locally conformal parallel \spin{7} structure on $S^7_\text{sq}\times S^1$, that is not
conformally flat.
\end{rem}

We restrict now our attention to compact locally conformal parallel flat $\spin{7}$ manifolds covered by
$\R^8-\{0\}$. The \spin{7} invariant 4-form $\sform$ on $\R^8$ is given, in terms of octonion multiplication, by
$\sform(x,y,z,v)=\langle(x\cdot y) \cdot z,v\rangle$.
Then the 4-form
\[\tilde{\sform}(x,y,z,w)=\frac{\sform(x,y,z,w)}{(|x|^2+|y|^2+|z|^2+|w|^2)^2}\]
descends to a
locally conformal parallel flat $\spin{7}$ structure on $S^7\times S^1$. 
For any finite $\Gamma\subset\spin{7}$ that
acts freely on
$S^7$ we get an example of locally conformal parallel
$\spin{7}$ structure on $(S^7/\Gamma)\times S^1$ as
well as a new nearly parallel $G_2$ structure (respectively
Sasaki-Einstein) on the factor $S^7/\Gamma$.

A class of finite subgroups of \spin{7} can be described by recalling that 
$\spin{7}$ is generated by the right multiplication on the octonions $\Ca$ by
elements of $S^6\subset \im{\Ca}$. Let $1\le m\le 7$ and let $V_m(\R^7)$ be the Stiefel manifold of
orthonormal $m$-frames in $\R^7$. Any $\sigma\in V_m(\R^7)$ gives rise to mutually orthogonal
$\sigma_1,\dots,\sigma_m\in S^6$ generating a finite
subgroup $G_{\sigma(m)}\subset\spin{7}$ of order $2^{m+1}$. This is due to
the identity $(o\bar y)x=-(o\bar x)y$
that holds for any $o\in\Ca$ and any orthogonal octonions $x,y$.
For example, the choice of $\sigma(4)=\{i,j,k,l\}$ gives rise to the group
$G_{\sigma(4)}$:
\begin{gather*}
o\rightarrow \pm o,\quad o\rightarrow \pm oi,\quad o\rightarrow \pm oj,\quad
o\rightarrow \pm ok, \quad o\rightarrow \pm ol,\quad o\rightarrow \pm (oi)j,\quad
o\rightarrow \pm (oi)k, \\ \nonumber
o\rightarrow \pm (oi)l,\quad
o\rightarrow \pm (oj)k,\quad o\rightarrow \pm (oj)l,\quad
o\rightarrow \pm (ok)l, o\rightarrow \pm ((oi)j)k,\\ \nonumber
o\rightarrow \pm ((oi)k)l, \quad o\rightarrow \pm ((oi)j)l,\quad
o\rightarrow \pm ((oj)k)l,\quad o\rightarrow \pm (((oi)j)k)l.
\end{gather*}

Other examples include the 
finite subgroup of $\symp{2}\subset\spin{7}$ obtained by
looking at $\Ca=\Ha\oplus\Ha$. A simple example is given by the
group:
\begin{equation*}
o=(q,q')\rightarrow (\pm (q,q'),\pm (qi,q'i),\pm (qj,q'j),\pm(qk,q'k)).
\end{equation*}

\begin{rem}
Note that most of these finite groups do not act freely on $S^7$
and it would be interesting to recognize which of the quotients
are actually smooth and can appear as spherical space forms, thus
fitting in J.~Wolf classification \cite{WolSCC,WolISS}. For any
such space form $S^7/\Gamma$ the product manifold
$(S^7/\Gamma)\times S^1$ will carry a locally conformal parallel
flat \spin{7} structure. Moreover, all compact locally conformal
parallel flat \spin{7} manifold could be
obtained in this way selecting from the J.~Wolf classification
\cite{WolSCC,WolISS} of finite subgroups of $SO(8)$ those which
are finite subgroups of \spin{7}.
\end{rem}

Another way of producing examples is to look at the locally conformal
parallel \spin{7} and non locally conformal K\"ahler Ricci flat 
$S_\text{sq}^7\times S^1$. Since the isometry group of
the squashed sphere $S_\text{sq}^7$ is $\symp{2}\cdot\symp{1}$, any finite 
$\Gamma\subset\symp{2}\cdot\symp{1}$ acting freely gives rise to such a structure on
$(S_\text{sq}^7/\Gamma)\times S^1$. These are pure examples in the sense that they are
not locally conformal K\"ahler Ricci flat. More generally, any 3-Sasakian structure
in dimension 7 gives rise to a `squashed' nearly parallel $\gtwo$ structure
\cite{FKMNPG}, that admits exactly one Killing spinor and similar arguments apply.


\medskip
\noindent Stefan Ivanov, University of Sofia \\   
Faculty of Mathematics and Informatics, Blvd. James Bourchier $5$,
1164 Sofia, Bulgaria, ivanovsp@fmi.uni-sofia.bg

\medskip
\noindent Maurizio Parton, Universit\`a di Chieti--Pescara\\
Dipartimento di Scienze, viale Pindaro $87$, I-65127 Pescara,
Italy,\\ parton@sci.unich.it

\medskip
\noindent Paolo Piccinni, Universit\`a degli Studi di Roma ``La
Sapienza''\\ Dipartimento di Matematica ``Guido Castelnuovo'', P.le
Aldo Moro 2, 00185, Roma, Italy,\\ piccinni@mat.uniroma1.it

\end{document}